\documentclass[onecolumn,11pt]{article}
\usepackage[top=1in, bottom=1in, left=1in, right=1in]{geometry}
\usepackage{amsmath,amsfonts,amscd,amssymb}
\usepackage{graphicx}
\graphicspath{{figures/}}
\usepackage{url}
\usepackage{xcolor}
\usepackage{caption}
\usepackage{mathtools}
\mathtoolsset{showonlyrefs}
\usepackage{setspace}
\usepackage{tikz}
\usetikzlibrary{shapes,arrows,calc}
\usepackage[numbers,sort&compress]{natbib}

\usepackage{listings}
\usepackage{xcolor}

\definecolor{maroon}{cmyk}{0, 0.87, 0.68, 0.32}
\definecolor{halfgray}{gray}{0.55}
\definecolor{ipython_frame}{RGB}{207, 207, 207}
\definecolor{ipython_bg}{RGB}{247, 247, 247}
\definecolor{ipython_red}{RGB}{186, 33, 33}
\definecolor{ipython_green}{RGB}{0, 128, 0}
\definecolor{ipython_cyan}{RGB}{64, 128, 128}
\definecolor{ipython_purple}{RGB}{170, 34, 255}

\lstdefinelanguage{iPython}{
    morekeywords={access,and,break,class,continue,def,del,elif,else,except,exec,finally,for,from,global,if,import,in,is,lambda,not,or,pass,print,raise,return,try,while,True,False,as},%
    morekeywords=[2]{abs,all,any,basestring,bin,bool,bytearray,callable,chr,classmethod,cmp,compile,complex,delattr,dict,dir,divmod,enumerate,eval,execfile,file,filter,float,format,frozenset,getattr,globals,hasattr,hash,help,hex,id,input,int,isinstance,issubclass,iter,len,list,locals,long,map,max,memoryview,min,next,object,oct,open,ord,pow,property,range,raw_input,reduce,reload,repr,reversed,round,set,setattr,slice,sorted,staticmethod,str,sum,super,tuple,type,unichr,unicode,vars,xrange,zip,apply,buffer,coerce,intern},%
    sensitive=true,%
    morecomment=[l]\#,%
    morestring=[b]',%
    morestring=[b]",%
    morestring=[s]{'''}{'''},
    morestring=[s]{"""}{"""},
    morestring=[s]{r'}{'},
    morestring=[s]{r"}{"},%
    morestring=[s]{r'''}{'''},%
    morestring=[s]{r"""}{"""},%
    morestring=[s]{u'}{'},
    morestring=[s]{u"}{"},%
    morestring=[s]{u'''}{'''},%
    morestring=[s]{u"""}{"""},%
    literate=
    {á}{{\'a}}1 {é}{{\'e}}1 {í}{{\'i}}1 {ó}{{\'o}}1 {ú}{{\'u}}1
    {Á}{{\'A}}1 {É}{{\'E}}1 {Í}{{\'I}}1 {Ó}{{\'O}}1 {Ú}{{\'U}}1
    {à}{{\`a}}1 {è}{{\`e}}1 {ì}{{\`i}}1 {ò}{{\`o}}1 {ù}{{\`u}}1
    {À}{{\`A}}1 {È}{{\'E}}1 {Ì}{{\`I}}1 {Ò}{{\`O}}1 {Ù}{{\`U}}1
    {ä}{{\"a}}1 {ë}{{\"e}}1 {ï}{{\"i}}1 {ö}{{\"o}}1 {ü}{{\"u}}1
    {Ä}{{\"A}}1 {Ë}{{\"E}}1 {Ï}{{\"I}}1 {Ö}{{\"O}}1 {Ü}{{\"U}}1
    {â}{{\^a}}1 {ê}{{\^e}}1 {î}{{\^i}}1 {ô}{{\^o}}1 {û}{{\^u}}1
    {Â}{{\^A}}1 {Ê}{{\^E}}1 {Î}{{\^I}}1 {Ô}{{\^O}}1 {Û}{{\^U}}1
    {œ}{{\oe}}1 {Œ}{{\OE}}1 {æ}{{\ae}}1 {Æ}{{\AE}}1 {ß}{{\ss}}1
    {ç}{{\c c}}1 {Ç}{{\c C}}1 {ø}{{\o}}1 {å}{{\r a}}1 {Å}{{\r A}}1
    {€}{{\EUR}}1 {£}{{\pounds}}1,
    literate=
    *{-}{{{\color{ipython_purple}-}}}1
     {?}{{{\color{ipython_purple}?}}}1
     {^}{{{\color{ipython_purple}\^{}}}}1
     {=}{{{\color{ipython_purple}=}}}1
     {+}{{{\color{ipython_purple}+}}}1
     {*}{{{\color{ipython_purple}$^\ast$}}}1
     {/}{{{\color{ipython_purple}/}}}1
     {>}{{{\color{ipython_purple}$>$}}}1
     {<}{{{\color{ipython_purple}$<$}}}1
     {@}{{{\color{ipython_purple}@}}}1
     {+=}{{{+=}}}1
     {-=}{{{-=}}}1
     {*=}{{{$^\ast$=}}}1
     {/=}{{{/=}}}1,
    identifierstyle=\color{black}\ttfamily,
    commentstyle=\color{ipython_cyan}\ttfamily,
    stringstyle=\color{ipython_red}\ttfamily,
    keepspaces=true,
    showspaces=false,
    showstringspaces=false,
    rulecolor=\color{ipython_frame},
    backgroundcolor=\color{ipython_bg},
    basicstyle=\normalsize,
    keywordstyle=\color{ipython_green}\ttfamily,
}

\definecolor{mygray}{gray}{0.95}
\definecolor{codegreen}{rgb}{0,0.6,0}
\definecolor{codegray}{rgb}{0.5,0.5,0.5}
\definecolor{codepurple}{rgb}{0.58,0,0.82}
\definecolor{backcolour}{rgb}{0.95,0.95,0.92}

\lstdefinestyle{mystyle}{
    backgroundcolor=\color{mygray},   
    commentstyle=\color{codegreen},
    keywordstyle=\color{magenta},
    numberstyle=\tiny\color{codegray},
    stringstyle=\color{codepurple},
    basicstyle=\ttfamily\footnotesize,
    breakatwhitespace=false,         
    breaklines=true,                 
    captionpos=b,                    
    keepspaces=true,                 
    showspaces=false,                
    showstringspaces=false,
    showtabs=false,                  
    tabsize=2
}

\usepackage[bottom,flushmargin,hang,multiple]{footmisc}
\usepackage{lipsum}
\newcommand\blfootnote[1]{%
  \begingroup
  \renewcommand\thefootnote{}\footnote{#1}%
  \addtocounter{footnote}{-1}%
  \endgroup
}

\definecolor{header1}{cmyk}{0,0,0,1}

\newcommand{\Xv}{\mathbf{X}}

\newcommand{\Thetav}{\boldsymbol{\Theta}}

\newcommand{\Xiv}{\boldsymbol{\Xi}}

\tikzstyle{block} = [rectangle, draw, fill=blue!20, 
    text width=4.5em, text centered, rounded corners, minimum height=4em]
\tikzstyle{line} = [draw, -latex']

\title{\vspace{-.65in}{\fontsize{16}{16}\selectfont \textbf{PySINDy: A Python package for the Sparse Identification of Nonlinear Dynamics from Data}}\vspace{-.15in}}

\author{\normalsize{
Brian M. de Silva$^{1*}$, Kathleen Champion$^{1*}$, Markus Quade$^2$,}\\ \normalsize{Jean-Christophe Loiseau$^3$, J. Nathan Kutz$^1$, Steven L. Brunton$^4$}\\
\footnotesize{$^1$ Department of Applied Mathematics, University of Washington, Seattle, WA 98195, United States}\\
\footnotesize{$^2$ Ambrosys GmbH, Potsdam, Germany}\\
\footnotesize{$^3$ ParisTech, Paris, France}\\
\footnotesize{$^4$ Department of Mechanical Engineering, University of Washington, Seattle, WA 98195, United States\vspace{-.2in}}
}
\date{}
\begin{document}
\maketitle

\blfootnote{$^*$ Corresponding authors (bdesilva@uw.edu and kpchamp@uw.edu).}
\vspace{-.2in}
\begin{abstract}
\texttt{PySINDy} is a Python package for the discovery of governing dynamical systems models from data.
In particular, \texttt{PySINDy} provides tools for applying the sparse identification of nonlinear dynamics (SINDy)~\cite{brunton2016pnas} approach to model discovery. In this work we provide a brief description of the mathematical underpinnings of SINDy, an overview and demonstration of the features implemented in \texttt{PySINDy} (with code examples),
practical advice for users, and a list of potential extensions to \texttt{PySINDy}.  Software is available at \url{https://github.com/dynamicslab/pysindy}.

\vspace{0.05in}
\vspace{-.15in}
\end{abstract}

\section{Introduction}

Scientists have long quantified empirical observations by developing mathematical models that characterize the observations, have some measure of interpretability, and are capable of making predictions.
Dynamical systems models, in particular, have been widely used to study, explain, and predict behavior in a diversity of application areas, with examples ranging from Newton's laws of classical mechanics to the Michaelis-Menten kinetics for modeling enzyme kinetics.
While governing laws and equations have traditionally been derived from first principles and expert knowledge, the current growth of available measurement data and the resulting emphasis on data-driven modeling motivates algorithmic and reproducible approaches for automated model discovery.

A number of such approaches have been developed in recent years~\cite{Brunton2019book}, including linear methods~\cite{Nelles2013book,ljung2010arc}, dynamic mode decomposition (DMD)~\cite{schmid2010dynamic,Kutz2016book} and Koopman theory more generally~\cite{Budivsic2012chaos,Mezic2013arfm,Williams2015jnls,klus2017data,Li2017chaos,Brunton2017natcomm}, nonlinear autoregressive algorithms~\cite{Akaike1969annals,Billings2013book}, neural networks~\cite{long2017pde,yang2018physics,Wehmeyer2018jcp,Mardt2018natcomm,vlachas2018data,pathak2018model,lu2019deepxde,Raissi2019jcp,Champion2019pnas,raissi2020science}, Gaussian process regression~\cite{Raissi2017arxiva,Raissi2017arxiv}, operator inference and reduced-order modeling~\cite{Benner2015siamreview,peherstorfer2016data,qian2020lift}, nonlinear Laplacian spectral analysis~\cite{Giannakis2012pnas}, diffusion maps~\cite{Yair2017pnas}, genetic programming~\cite{bongard_automated_2007,schmidt_distilling_2009,Daniels2015naturecomm}, and sparse regression~\cite{brunton2016pnas,Rudy2017sciadv}.
Maximizing the impact of these model discovery methods requires tools to make them widely accessible to scientists across domains and at various levels of mathematical expertise.

\texttt{PySINDy} is a Python package for the discovery of governing dynamical systems models from data. In particular, \texttt{PySINDy} provides tools for applying the sparse identification of nonlinear dynamics (SINDy) approach to model discovery~\cite{brunton2016pnas}. SINDy poses model discovery as a sparse regression problem, where relevant terms in the dynamics are selected from a library of candidate functions, many of them motivated by our deep historical knowledge of diverse physics models. 
This approach results in interpretable models, and it has been widely applied~\cite{Sorokina2016oe,Loiseau2017jfm,Dam2017pf,Loiseau2018jfm,Hoffmann2018arxiv,Loiseau2019data,el2018sparse,narasingam2018data,de2019discovery,Thaler2019jcp,lai2019sparse,Deng2020JFM,schmelzer2020discovery,pan2020sparsity,beetham2020formulating}  and extended~\cite{Mangan2016ieee,Loiseau2017jfm,Schaeffer2017prsa,Mangan2017prsa,Rudy2017sciadv,tran2017exact,Schaeffer2017pre,schaeffer2018extracting,wu2018numerical,boninsegna2018sparse,Kaiser2018prsa,mangan2019model,Gelss2019mindy,goessmann2020tensor,Reinbold2020pre} using different sparse optimization algorithms and library functions.

The \texttt{PySINDy} package is aimed at researchers and practitioners alike, enabling anyone with access to measurement data to engage in scientific model discovery. The package is designed to be accessible to inexperienced users, adhering to \texttt{scikit-learn} standards, while also including customizable options for more advanced users.  
A number of popular SINDy variants are implemented, but \texttt{PySINDy} is also designed to enable further extensions for research and experimentation.  

\section{Background}\label{sec:background}

\texttt{PySINDy} provides an implementation of the SINDy method to discover governing dynamical systems models of the form 
\begin{equation}\label{eq:dynamical_system}
  \frac{d}{dt}\mathbf{x}(t) = \mathbf{f}(\mathbf{x}(t)). 
\end{equation}
Given data in the form of state measurements $\mathbf{x}(t) \in \mathbb{R}^n$, SINDy identifies a model for the dynamics, given by the function $\mathbf{f}$, which describes how the state of the system evolves in time.  In particular, SINDy sparsely approximates the dynamics in a library of candidate basis functions $\boldsymbol{\theta}(\mathbf{x}) = [\theta_1(\mathbf{x}),\theta_2(\mathbf{x}),\dots,\theta_\ell(\mathbf{x})]$, so that 
\begin{equation}
\mathbf{f}(\mathbf{x})\approx \sum_{k=1}^{\ell}\theta_k(\mathbf{x})\xi_k.
\end{equation}
The majority of coefficients $\xi_k$ are zero, and nonzero entries identify active terms in the dynamics.  

To pose SINDy as a regression problem, time-series measurements of $\mathbf{x}$ and their time derivatives $\dot{\mathbf{x}}$ are arranged into matrices
\begin{equation*}
 \mathbf{X} = \left(\begin{array}{cccc}
   x_1(t_1) & x_2(t_1) & \cdots & x_n(t_1) \\
   x_1(t_2) & x_2(t_2) & \cdots & x_n(t_2) \\
   \vdots & \vdots & \ddots & \vdots \\
   x_1(t_m) & x_2(t_m) & \cdots & x_n(t_m) \\
 \end{array}\right), \quad
 \dot{\mathbf{X}} = \left(\begin{array}{cccc}
   \dot{x}_1(t_1) & \dot{x}_2(t_1) & \cdots & \dot{x}_n(t_1) \\
   \dot{x}_1(t_2) & \dot{x}_2(t_2) & \cdots & \dot{x}_n(t_2) \\
   \vdots & \vdots & \ddots & \vdots \\
   \dot{x}_1(t_m) & \dot{x}_2(t_m) & \cdots & \dot{x}_n(t_m) \\
 \end{array}\right).
\end{equation*}
The derivatives can be approximated numerically or measured directly.
The library functions are evaluated on the data, resulting in ${\boldsymbol{\Theta}(\mathbf{X}) = [\theta_1(\mathbf{X}),\theta_2(\mathbf{X}),\dots,\theta_\ell(\mathbf{X})]}$.  Sparse regression is then performed to approximately solve
\begin{equation}
  \dot{\mathbf{X}} \approx \boldsymbol{\Theta}(\mathbf{X}) \boldsymbol{\Xi}, \label{eq:sindy_regression}
\end{equation}
where $\boldsymbol{\Xi}$ is a set of coefficients that determines the active terms in $\mathbf{f}$. While the original SINDy formulation solves~\eqref{eq:sindy_regression} using a sequentially thresholded least squares algorithm~\cite{brunton2016pnas,zhang2019convergence}, this problem can be solved using any sparse regression algorithm, such as lasso~\cite{TibshiraniLasso},  sparse relaxed regularized regression~(SR3)~\cite{zheng2018ieee,champion2019}, stepwise sparse regression~(SSR)~\cite{boninsegna2018sparse}, or Bayesian methods~\cite{Guang2018,Pan2016BayesianSINDy,niven2020bayesian}.

SINDy has been widely applied for model identification in applications such as chemical reaction dynamics~\cite{Hoffmann2018arxiv}, nonlinear optics~\cite{Sorokina2016oe}, fluid dynamics~\cite{Loiseau2017jfm,Loiseau2018jfm,Loiseau2019data,el2018sparse,Deng2020JFM} and turbulence modeling~\cite{schmelzer2020discovery,beetham2020formulating}, plasma convection~\cite{Dam2017pf}, numerical algorithms~\cite{Thaler2019jcp}, and structural modeling~\cite{lai2019sparse}, among others~\cite{narasingam2018data,de2019discovery,pan2020sparsity}.
It has also been extended to handle more complex modeling scenarios such as partial differential equations~\cite{Schaeffer2017prsa,Rudy2017sciadv}, systems with inputs or control~\cite{Kaiser2018prsa}, systems with implicit dynamics~\cite{Mangan2016ieee}, hybrid systems~\cite{mangan2019model}, to enforce physical constraints~\cite{Loiseau2017jfm}, to incorporate information theory~\cite{Mangan2017prsa}, to identify models from corrupt or limited data~\cite{tran2017exact,schaeffer2018extracting} and ensembles of initial conditions~\cite{wu2018numerical}, and extending the formulation to include integral terms~\cite{Schaeffer2017pre,Reinbold2020pre}, tensor representations~\cite{Gelss2019mindy,goessmann2020tensor}, and stochastic forcing~\cite{boninsegna2018sparse}.
However, there is not a definitive standard implementation or package for SINDy.
Versions of SINDy have been implemented within larger projects such as \texttt{sparsereg}~\cite{markus_quade_sparsereg}, but no specific implementation has emerged as the most widely adopted and most versions implement a limited set of features.
Researchers have thus typically written their own implementations, resulting in duplicated effort and a lack of standardization. This not only makes it more difficult to apply SINDy to scientific data sets, but also makes it more challenging to benchmark extensions of the method and makes such extensions less accessible to end users. This motivates the creation of a dedicated package for SINDy. The \texttt{PySINDy} package provides a central codebase where many of the basic SINDy features are implemented, allowing for easy use and standardization. In addition, it is straightforward for users to extend \texttt{PySINDy} so that new developments are available to the wider community.

\section{Features}
The core object in the \texttt{PySINDy} package is the \texttt{SINDy} model class, which is implemented as a \texttt{scikit-learn} estimator. This design choice was made to ensure that the package is simple to use for a wide user base, as many potential users will be familiar with \texttt{scikit-learn}. It also expresses the \texttt{SINDy} model object at the appropriate level of abstraction so that users can embed it into more sophisticated pipelines in \texttt{scikit-learn}, such as for parameter tuning and model selection.

Our \texttt{PySINDy} implementation involves three major steps, resulting in three modeling decisions:
\begin{enumerate}
  \item The \textit{numerical differentiation} scheme used to approximate $\dot{\Xv}$ from $\Xv$;
  \item The candidate functions constituting the \textit{feature library} $\Thetav$;
  \item The \textit{sparse regression} algorithm that is applied to solve \eqref{eq:sindy_regression} to find $\Xiv$.
\end{enumerate}

\begin{figure}
  \centering
  \begin{tikzpicture}[node distance=3cm]
    \node [block, line width=1.5pt, fill=black!20] (xdot) {\Large$\dot{\Xv}$};
    \node [right of=xdot, line width=1.5pt, node distance=2cm] (approx) {\Large$\approx$};
    \node [block, line width=1.5pt, right of=approx, node distance=2cm] (theta) {\Large$\Thetav(\Xv)$};
    \node [block, line width=1.5pt, right of=theta, fill=red!20, node distance=2.2cm] (xi) {\Large$\Xiv$};

    \node [block, line width=1.5pt, below of=theta, text width=14em, node distance=2.5cm] (library) {\texttt{pysindy.feature\_library}};
    \node [block, line width=1.5pt, left of=library, fill=black!20, text width=14em, node distance=6.1cm] (diff) {\texttt{pysindy.differentiation}};
    \node [block, line width=1.5pt, right of=library, fill=red!20, text width=11em, node distance=5.5cm] (opt) {\texttt{pysindy.optimizers}};

    \path [line, line width=1.5pt] (xdot) -| (diff);
    \path [line, line width=1.5pt] (theta) -- (library);
    \path [line, line width=1.5pt] (xi) -| (opt);

    \node [below of=diff, node distance=1.2cm, line width=1.5pt] {Numerical differentiation};
    \node [below of=library, node distance=1.2cm, line width=1.5pt] {Form feature library};
    \node [below of=opt, node distance=1.2cm, line width=1.5pt] {Sparse regression};

  \end{tikzpicture}
  \caption{Correspondence between the sparse regression problem solved by SINDy and the submodules of PySINDy.}
  \label{fig:package-structure}
\end{figure}
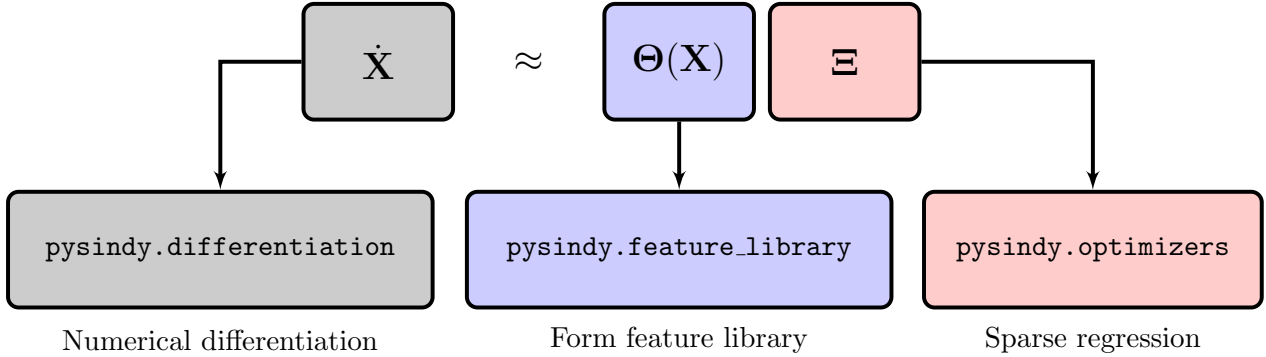

The core \texttt{SINDy} object was designed to incorporate these three components in as modular a manner as possible, having one attribute corresponding to each: \texttt{SINDy.differentiation\_method} for numerical differentiation, \texttt{SINDy.feature\_library} for the formation of the candidate function library, and \texttt{SINDy.optimizer} for the sparse regressor. 
\texttt{PySINDy} provides standard options for each step, while making it easy to replace any of these steps with more sophisticated ``third-party'' algorithms. In particular, at the time of writing, we have implemented the following methods:
\begin{itemize}
  \item Numerical differentiation (for computing $\dot{\Xv}$ from $\Xv$)
  \begin{itemize}
    \item Finite difference: \texttt{FiniteDifference}
    \item Smoothed finite difference: \texttt{SmoothedFiniteDifference}
  \end{itemize}
  \item Feature library (for constructing $\Thetav$)
  \begin{itemize}
    \item Multivariate polynomials: \texttt{PolynomialLibrary}
    \item Fourier modes (i.e. trigonometric functions): \texttt{FourierLibrary}
    \item Custom library (defined by user-supplied functions): \texttt{CustomLibrary}
    \item Identity library (in case users want to compute $\Thetav$ themselves): \texttt{IdentityLibrary}
  \end{itemize}
  \item Optimizer (for performing sparse regression)
  \begin{itemize}
    \item Sequentially thresholded least-squares~\cite{brunton2016pnas,zhang2019convergence}: \texttt{STLSQ}
    \item Sparse relaxed regularized regression (SR3)~\cite{zheng2018ieee}: \texttt{SR3}
  \end{itemize}
\end{itemize}

\section{Examples}
The \texttt{PySINDy} GitHub page\footnote{\url{https://github.com/dynamicslab/pysindy}} includes tutorials in the form of Jupyter notebooks.
These tutorials demonstrate the usage of various features of the package and reproduce the examples from the original SINDy paper~\cite{brunton2016pnas}.
Throughout this section we will use the Lorenz equations \eqref{eq:lorenz} as the dynamical system to illustrate the \texttt{PySINDy} package:
\begin{equation}\label{eq:lorenz}
  \left\{
  \begin{aligned}
    \dot x &= -10x - 10y \\
    \dot y &= x(28 - z) - y \\
    \dot z &= xy - \frac83 z
  \end{aligned}
  \right.
\end{equation}

\begin{figure}
  \centering
  \includegraphics[width=.5\textwidth]{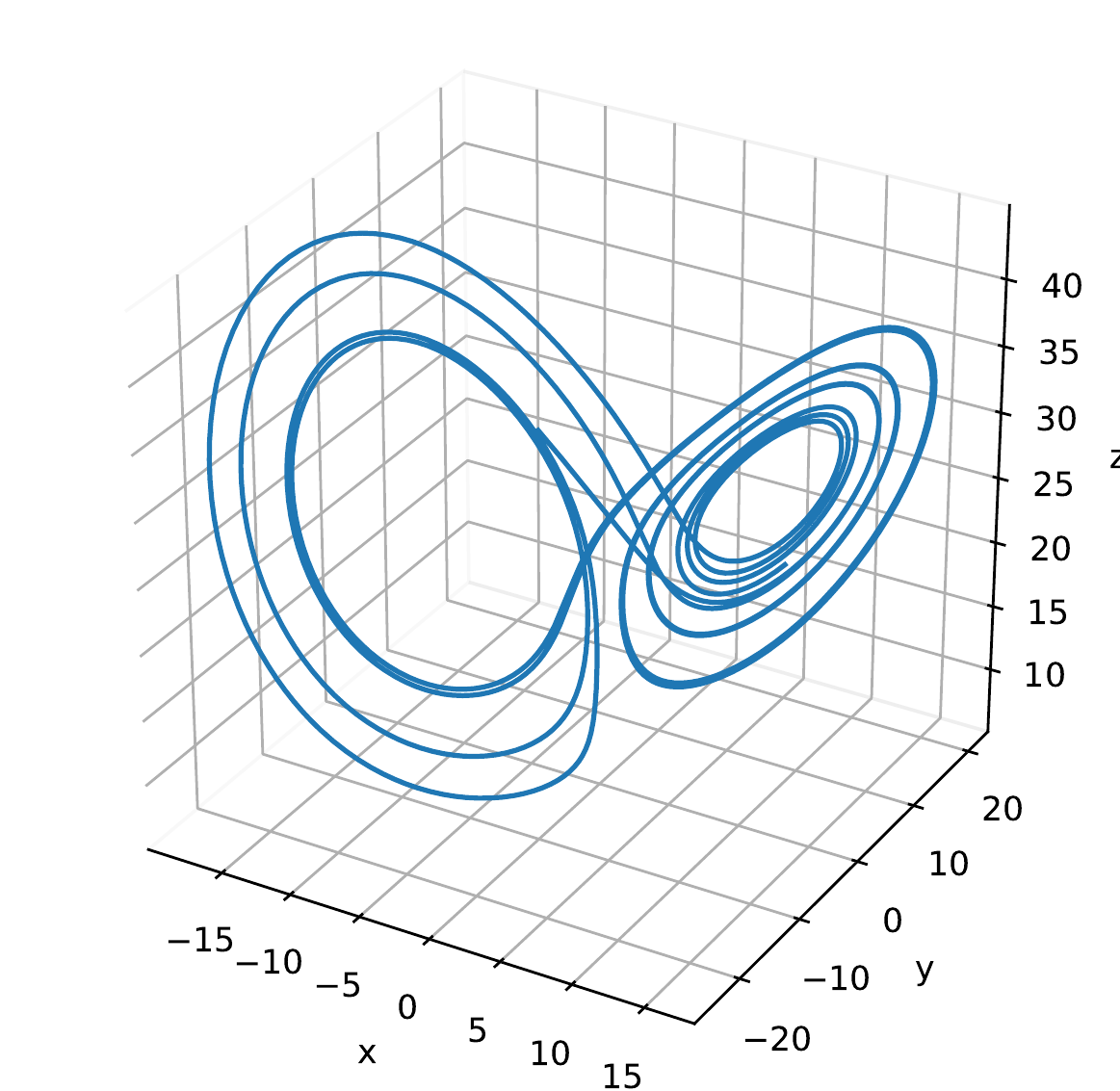}
  \caption{Measurement data simulated using the Lorenz equations \eqref{eq:lorenz}.} 
  \label{fig:lorenz_data}
\end{figure}

In Python, the right-hand side of \eqref{eq:lorenz} can be expressed as follows:
\begin{lstlisting}[language=iPython]
def lorenz(x, t):
    return [
        10 * (x[1] - x[0]),
        x[0] * (28 - x[2]) - x[1],
        x[0] * x[1] - (8 / 3) * x[2]
    ]
\end{lstlisting}

To construct training data to feed into a SINDy model, we integrate \eqref{eq:lorenz} with:
\begin{lstlisting}[language=iPython]
import numpy as np
from scipy.integrate import odeint

dt = 0.002
t = np.arange(0, 10, dt)
x0 = [-8, 8, 27]
X = odeint(lorenz, x0, t)
\end{lstlisting}

We plot \texttt{X} in Figure \ref{fig:lorenz_data}. It is important to note that each column of \texttt{X} corresponds to a variable and each row to a point in time. All \texttt{PySINDy} objects that handle data assume the data is structured this way.

\subsection{Basic usage}
  The \texttt{pysindy} package is built around the \texttt{SINDy} class, which encapsulates all the steps necessary to learn a dynamical system with SINDy. 
  To create a SINDy object, fit it to the data (i.e. to infer a dynamical system from the data), and print the resulting model, we invoke the \texttt{SINDy} constructor, the \texttt{fit} method, and custom \texttt{print} functions
  \begin{lstlisting}[language=iPython]
model = ps.SINDy()
model.fit(X, t=dt)
model.print()
\end{lstlisting}
  
 \noindent which generates the following output
  \begin{verbatim}
  x0' = -9.999 x0 + 9.999 x1
  x1' = 27.992 x0 + -0.999 x1 + -1.000 x0 x2
  x2' = -2.666 x2 + 1.000 x0 x1
  \end{verbatim}

  Once the SINDy object has been fit we can feed in new data and use the learned model to predict the derivatives for each measurement (recall that measurements correspond to rows).
\begin{lstlisting}[language=iPython]
t_test = np.arange(0, 15, dt)
x0_test = np.array([8, 7, 15])
X_test = odeint(lorenz, x0_test, t_test)

X_dot_test_computed = model.differentiate(X_test, t=dt)
X_dot_test_predicted = model.predict(X_test)
\end{lstlisting}
  
\begin{figure}
    \centering
    \includegraphics[width=\textwidth]{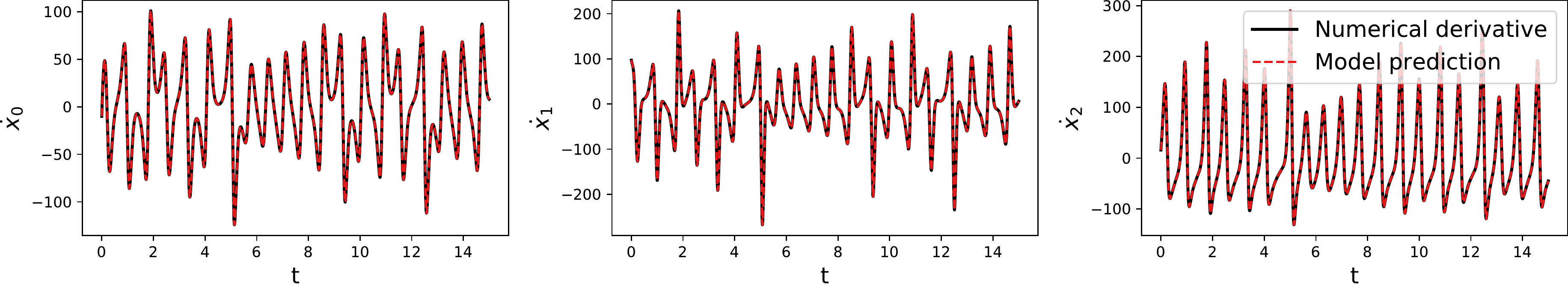}
    \caption{Derivatives of variables from the Lorenz equation via numerical differentiation and using a learned SINDy model.} 
    \label{fig:derivatives}
  \end{figure}

  The call \texttt{model.differentiate(X\_test, t=dt)} applies the numerical differentiation method in the \texttt{SINDy} model  to \texttt{X\_test} with time steps of length \texttt{dt}. In Figure \ref{fig:derivatives} we plot each dimension of \texttt{X\_dot\_test\_computed} and \texttt{X\_dot\_test\_predicted}.

  Rather than predicting derivatives, we will often be interested in using our model to evolve initial conditions forward in time using the learned model. The \texttt{simulate} function does just that.
  \begin{lstlisting}[language=iPython]
X_test_sim = model.simulate(x0_test, t_test)
\end{lstlisting}

\begin{figure}
    \centering
    \includegraphics[width=\textwidth]{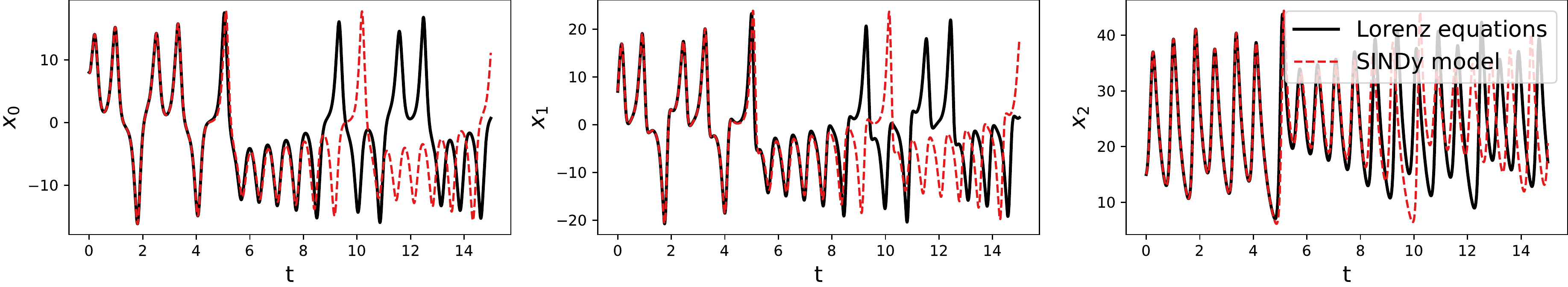}
    \caption{Two trajectories starting at the same position evolved forward in time with the exact Lorenz equations (black, solid) and the learned SINDy model (red, dashed).} 
    \label{fig:lorenz-simulated}
  \end{figure}

  Figure \ref{fig:lorenz-simulated} shows the simulated trajectory plotted against the true trajectory \texttt{X\_test}. The trajectories agree initially, but they eventually diverge due to the chaotic nature of the Lorenz equations.

\subsection{Custom features}
Thus far we have relied on the default options of the \texttt{SINDy} object, but \texttt{PySINDy} comes equipped with multiple alternative built-in methods for differentiation, library building, and optimization. These options are selected by passing corresponding \texttt{PySINDy} objects to the \texttt{SINDy} constructor via the \texttt{differentiation\_method}, \texttt{feature\_library}, and \texttt{optimizer} arguments, respectively. Parameters for the differentiation, library, and optimization algorithms are supplied to the corresponding objects' constructors rather than directly to the \texttt{SINDy} object. We demonstrate the syntax with the following example.
\begin{lstlisting}[language=iPython]
differentiation_method = ps.FiniteDifference(order=1)
feature_library = ps.PolynomialLibrary(degree=3, include_bias=False)
optimizer = ps.SR3(threshold=0.1, nu=1, tol=1e-6)

model = ps.SINDy(
    differentiation_method=differentiation_method,
    feature_library=feature_library,
    optimizer=optimizer,
    feature_names=["x", "y", "z"]
)

model.fit(X, t=dt)
model.print()
\end{lstlisting}

\noindent which prints
\begin{verbatim}
  x' = -10.021 x + 9.993 y
  y' = 28.431 x + -1.212 y + -1.008 x z
  z' = -2.675 z + 1.000 x y.
\end{verbatim}

A number of other built-in options are available. The official documentation\footnote{\url{https://pysindy.readthedocs.io/en/latest/index.html}} and examples\footnote{\url{https://github.com/dynamicslab/pysindy/tree/master/example}} provide an exhaustive list.

\newpage
\section{Practical tips}\label{sec:practical-tips}

  In this section we provide pragmatic advice for using \texttt{PySINDy} effectively. We discuss potential pitfalls and strategies for overcoming them. We also specify how to incorporate custom methods not implemented natively in \texttt{PySINDy}, where applicable. The information presented here is derived from a combination of experience and theoretical considerations.

  \subsection{Numerical differentiation}
    Numerical differentiation is one of the core components of the SINDy method. Derivatives of measurement variables provide the targets for the sparse regression problem~\eqref{eq:sindy_regression}. If care is not taken in computing these derivatives, the quality of the learned model is likely to suffer.

\begin{figure}[b]
      \centering
      \includegraphics[width=\textwidth]{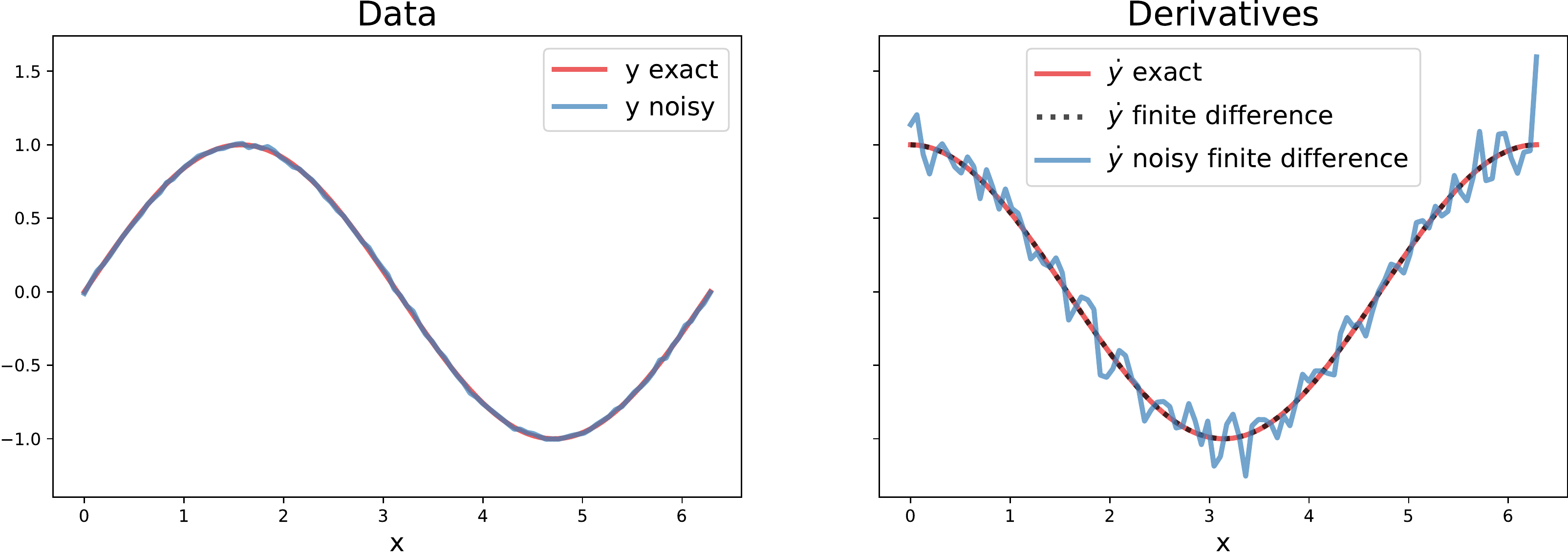}
      \caption{A toy example illustrating the effect of noise on derivatives computed with a second order finite difference method. Left: The data to be differentiated; $y=\sin(x)$ with and without a small amount of additive noise (normally distributed with mean 0 and standard deviation $0.01$). Right: Derivatives of the data; the exact derivative $\cos(x)$ (blue), the finite difference derivative of the exact data (black, dashed), and the finite difference derivative of the noisy data.}
      \label{fig:noisy_differentiation}
    \end{figure}

    By default, a second order finite difference method is used to differentiate input data. Finite difference methods tend to amplify noise in data. If the data are smooth (at least twice differentiable), then finite difference methods give accurate derivative approximations. When the data are noisy, they give derivative estimates with \textit{more} noise than the original data. Figure \ref{fig:noisy_differentiation} visualizes the impact of noise on numerical derivatives. Note that even a small amount of noise in the data can produce noticeable degradation in the quality of the numerical derivative.

    One way to mitigate the effects of noise is to smooth the measurements before computing derivatives. The \texttt{SmoothedFiniteDifference} method can be used for this purpose.
    A numerical differentiation scheme with total variation regularization has also been proposed~\cite{chartrand_numerical_2011} and recommended for use in SINDy~\cite{brunton2016pnas}.

    Users wishing to employ their own numerical differentiation schemes have two ways of doing so. Derivatives of input measurements can be computed externally with the method of choice and then passed directly into the \texttt{SINDy.fit} method via the \texttt{x\_dot} keyword argument. Alternatively, users can implement their own differentiation methods and pass them into the \texttt{SINDy} constructor using the \texttt{differentiation\_method} argument. In this case, the supplied class need only have implemented a \texttt{\_\_call\_\_} method taking two arguments, \texttt{x} and \texttt{t}.
  
  \subsection{Library selection}
    The SINDy method assumes dynamics can be represented as a \textit{sparse} linear combination of library functions. If this assumption is violated, the method is likely to exhibit poor performance. This issue tends to manifest itself as numerous library terms being active, often with weights of vastly different magnitudes, still resulting in poor model error.

    Typically, prior knowledge of the system of interest and its dynamics should be used to make a judicious choice of basis functions. When such information is unavailable, the default class of library functions, polynomials, are a good place to start, as smooth functions have rapidly converging Taylor series. Brunton et al.~\cite{brunton2016pnas} showed that, equipped with a  polynomial library, SINDy can recover the first few terms of the (zero-centered) Taylor series of the true right-hand side function $\mathbf{f}(x)$. If one has reason to believe the dynamics can be sparsely represented in terms of Chebyshev polynomials rather than monomials, then the library should include Chebyshev polynomials.

    \texttt{PySINDy} includes the \texttt{CustomLibrary} and \texttt{IdentityLibrary} objects to allow for flexibility in the library functions. When the desired library consists of a set of functions that should be applied to each measurement variable (or pair, triplet, etc. of measurement variables) in turn, the \texttt{CustomLibrary} class should be used. The \texttt{IdentityLibrary} class is the most customizable, but transfers the work of computing library functions over to the user. It expects that all the features one wishes to include in the library have already been computed and are present in \texttt{X} before \texttt{SINDy.fit} is called, as it simply applies the identity map to each variable that is passed to it. 
    It is best suited for situations in which one has very specific instructions for how to apply library functions (e.g. if some of the functions should be applied to only some of the input variables).

    As terms are added to the library, the underlying sparse regression problem becomes increasingly ill-conditioned. Therefore it is recommended to start with a small library whose size is gradually expanded until the desired level of performance is achieved. 
   For example, a user may wish to start with a library of linear terms and then add quadratic and cubic terms as necessary to improve model performance.  
    For the best results, the strength of regularization applied should be increased in proportion to the size of the library to account for the worsening condition number of the resulting linear system.

    Users may also choose to implement library classes tailored to their applications. To do so one should have the new class inherit from our \texttt{BaseFeatureLibrary} class. See the documentation for guidance on which functions the new class is expected to implement.
 
  \subsection{Optimization}
    \texttt{PySINDy} uses various optimizers to solve the sparse regression problem. For a fixed differentiation method, set of inputs, and candidate library, there is still some variance in the dynamical system identified by SINDY, depending on which optimizer is employed.

    The default optimizer in \texttt{PySINDy} is the sequentially-thresholded least-squares algorithm (\texttt{STLSQ}). In addition to being the method originally proposed for use with SINDy, it involves a single, easily interpretable hyperparameter, and it exhibits good performance across a variety of problems.

    The sparse relaxed regularized regression (\texttt{SR3})~\cite{zheng2018ieee,champion2019} algorithm can be used when the results of \texttt{STLSQ} are unsatisfactory. It involves a few more hyperparameters that can  be tuned for improved accuracy. In particular, the \texttt{thresholder} parameter controls the type of regularization that is applied. For optimal results, one may find it useful to experiment with $L^0$, $L^1$, and clipped absolute deviation (CAD) regularization. The other hyperparameters can be tuned with cross-validation.

    Custom or third party sparse regression methods are also supported. Simply instantiate an instance of the custom object and pass it to the \texttt{SINDy} constructor using the \texttt{optimizer} keyword. Our implementation is compatible with any of the linear models from Scikit-learn (e.g. \texttt{RidgeRegression}, \texttt{Lasso}, and \texttt{ElasticNet}).
    See the documentation for a list of methods and attributes a custom optimizer is expected to implement. There you will also find an example where the Scikit-learn \texttt{Lasso} object is used to perform sparse regression.

  \subsection{Regularization}
    Regularization, in this context, is a technique for improving the conditioning of ill-posed problems. Without regularization, one often obtains highly unstable results, with learned parameter values differing substantially for slightly different inputs. SINDy seeks weights that express dynamics as a \textit{sparse} linear combination of library functions. When the columns of the measurement data or the library are statistically correlated, which is likely for  large libraries, the SINDy inverse problem can quickly become ill-posed. Though the sparsity constraint is a type of regularization itself, for many problems another form of regularization is needed for SINDy to learn a robust dynamical model.

    In some cases regularization can be interpreted as enforcing a prior distribution on the model parameters~\cite{bishop2006pattern}.
    Applying strong regularization biases the learned weights \textit{away} from the values that would allow them to best fit the data and \textit{toward} the values preferred by the prior distribution (e.g. $L^2$ regularization corresponds to a Gaussian prior).
    Therefore once a sparse set of nonzero coefficients is discovered, our methods apply an extra  ``unbiasing'' step where \textit{unregularized} least-squares is used to find the values of the identified nonzero coefficients.
    All of our built-in methods use regularization by default.

    Some general best practices regarding regularization follow. Most problems will benefit from some amount of regularization. Regularization strength should be increased as the size of the candidate right-hand side library grows. If warnings about ill-conditioned matrices are generated when \texttt{SINDy.fit} is called, more regularization may help. We also recommend setting \texttt{unbias} to \texttt{True} when invoking the \texttt{SINDy.fit} method, especially when large amounts of regularization are being applied. Cross-validation can be used to select appropriate regularization parameters for a given problem.

\section{Extensions}
  In this section we list potential extensions and enhancements to our SINDy implementation. We provide references for the improvements that are inspired by previously conducted research and the rationale behind the other potential changes.
  
  \begin{itemize}
    \item \textbf{Partial differential equations (PDEs):} While dynamical systems given by ordinary differential equations (ODEs) provide a flexible approach to modeling physical systems, many systems are inherently described by partial differential equations (PDEs), which are not immediately discoverable using \texttt{PySINDy}. 
    Multiple approaches for the data-driven discovery of PDEs have been proposed~\cite{Rudy2017sciadv,Schaeffer2017prsa} as extensions to the SINDy method, and these may be readily included within the \texttt{PySINDy} framework.
    \item \textbf{Identifying coordinates and latent variables:}  Many complex systems, such as fluid flows, are high-dimensional, yet exhibit low-dimensional patterns that may be exploited for modeling~\cite{Taira2017aiaa,Taira2020aiaaj,Brunton2020arfm}.  Identifying effective coordinate systems on which to build models is an important aspect of data-driven discovery.  Recently, SINDy has been embedded into an autoencoder framework~\cite{Champion2019pnas} to simultaneously identify effective coordinates and sparse dynamics.  Similarly, for many systems, it is impossible to measure the full state of the system, so that there are latent variables.  Time-delay coordinates have been useful for identifying sparse models from limited measurements~\cite{Brunton2017natcomm}.  Both of these are candidates for future extensions. 
    \item \textbf{Constraints:} SINDy has been extended to enforce \textit{physical constraints} during the sparse regression step~\cite{Loiseau2017jfm}. When working with physical systems with known conserved quantities, such as the conservation of energy in incompressible fluids~\cite{Loiseau2017jfm}, such a method allows one to automatically incorporate this prior information into the model discovery process. 
    \item \textbf{Integral formulation:} We previously discussed how measurement data with too much noise can disrupt the model discovery process, and we offered smoothing as one possible solution. Another is to work with an integral version of \eqref{eq:dynamical_system}, as proposed by Schaeffer and McCalla~\cite{Schaeffer2017pre} and extended to PDEs by Reinbold, Gurevich, and Grigoriev~\cite{Reinbold2020pre}. Where numerical differentiation tends to amplify noise, numerical integration tends to smooth it out. This formulation has been shown to improve the robustness of SINDy to noise.
    \item \textbf{Ensembles:} Ensembles are a proven class of methods in machine learning for variance reduction (improved model generalizability) at the cost of extra computation. Rather than training a single model, multiple high-variance models are trained and their predictions are averaged together. We think that ideas from ensemble learning could be adapted to improve the performance of SINDy models. Recent work has hinted at possible approaches~\cite{sachdeva2019pyuoi}.
    \item \textbf{Extended libraries:} Choosing the appropriate basis in which to represent dynamics is of critical importance for the successful application of SINDy. Although we currently provide methods allowing users the flexibility to input their own library functions, we aim to make the library construction process even easier by providing a common suite of tools for the creation and combination of sets of candidate functions. More basis functions could be supported natively, such as nonautonomous terms (those depending explicitly on the dependent variable, time). Taking this idea a step further, specific variables (columns of $\mathbf{X}$) which should not appear on the left-hand side of \eqref{eq:sindy_regression} could be identified by the user. This would enable SINDy to include inputs and control variables~\cite{Kaiser2018prsa}. Operations acting on one or more libraries could also be implemented. For example, combining libraries via union, intersection, composition, or tensor product could enable the expression of complicated nonlinear dynamics. The ability to apply libraries to only subsets of state variables could help cut down on the computational cost and improve the conditioning of the sparse regression problem solved within SINDy.
  \end{itemize}

\section{Acknowledgments}
\texttt{PySINDy} is a fork of \texttt{sparsereg}~\cite{markus_quade_sparsereg}. 
SLB acknowledges funding support from the Army Research Office (ARO W911NF-19-1-0045)  and the Air Force Office of Scientific Research (AFOSR FA9550- 18-1-0200). JNK acknowledges support from the Air Force Office of Scientific Research (AFOSR FA9550-17-1-0329). This material is based upon work supported by the National Science Foundation Graduate Research Fellowship under Grant Number DGE- 1256082.

\newpage
\begin{spacing}{.9}
  \small{
    \setlength{\bibsep}{6.5pt}
\bibliographystyle{IEEEtran}
    \bibliography{references}
  }
\end{spacing}

\end{document}